\documentclass[10pt, a4paper,twoside]{article}
\usepackage{enumerate}
\usepackage{color}
\usepackage{amsmath}

\usepackage{pb-diagram}
\usepackage[english]{babel}
\usepackage{amssymb,latexsym}
\usepackage{latexsym}

\def \bui#1#2{\mathrel{\mathop{\kern 0pt#1}\limits^{#2}}}
\newcommand{\R}{{\mathbb R}}

\newcommand{\lquot}[2]{\raisebox{-0.5ex}{$#2$}\backslash\!\raisebox{0.5ex}{$#1$}}

\newtheorem{example}{Examples}[section]
\newtheorem{thm}{Theorem}[section]
\newtheorem{lemma}[thm]{Lemma}
\newtheorem{prop}[thm]{Proposition}

\newtheorem{cor}[thm]{Corollary}
\newtheorem{remark}[thm]{Remark}
\newtheorem{remarks}[thm]{Remarks}
\newtheorem{definition}[thm]{Definition}
\newtheorem{notation}[thm]{Notation}
\newtheorem{exabout:ample}[thm]{Example}

\title{Rigidity results for spin manifolds with foliated boundary}
\author{Fida El Chami\footnote{Lebanese University, Faculty of Sciences II, Department of Mathematics, P.O. Box 90656 Fanar-Matn, Lebanon,
E-mail: \texttt{fchami@ul.edu.lb}},\, Nicolas Ginoux\footnote{IUT de Metz, D\'epartement informatique, Ile du Saulcy, CS 10628, 57045 Metz Cedex 01, France, E-mail: \texttt{nicolas.ginoux@univ-lorraine.fr}},\, Georges Habib\footnote{Lebanese University, Faculty of Sciences II, Department of Mathematics, P.O. Box 90656 Fanar-Matn, Lebanon,
E-mail: \texttt{ghabib@ul.edu.lb}},\, Roger Nakad \footnote{Notre Dame University-Louaiz\'e, Faculty of Natural and Applied Sciences, Department of Mathematics and Statistics, P.O. Box 72, Zouk Mikael, Lebanon, E-mail: \texttt{rnakad@ndu.edu.lb}}}

\begin{document}
\date{}
\maketitle
\begin{abstract}
\noindent In this paper, we consider a compact Riemannian manifold whose boundary is endowed with a Riemannian flow. Under a suitable curvature assumption depending on the O'Neill tensor of the flow, we prove that any solution of the basic Dirac equation is the restriction of a parallel spinor field defined on the whole manifold. As a consequence, we show that the flow is a local product. In particular, in the case where
solutions of the basic Dirac equation are given by basic Killing spinors, we characterize the geometry of the manifold and the flow.
\end{abstract}


\section{Introduction}

\noindent Over the last decades, spin geometry has become a useful
tool for the study of the geometry of manifolds and their
hypersurfaces. In this setting, many results have been obtained in
several papers \cite{B,HM,HM1,HM2,HM3, H-M} relating the extrinsic
aspect of an hypersurface $M$ embedded (or immersed) into a
Riemannian manifold $N$, to the intrinsic one.

\noindent Using the spinorial Reilly formula \cite[Remark 1]{HM3},
O.~Hijazi, S.~Montiel and X.~Zhang proved that on the compact
boundary $M^{n+1}$ of a  Riemannian spin manifold $N^{n+2}$ of
dimension $n + 2$ with non-negative scalar curvature, the first
positive eigenvalue $\lambda_1$ of the Dirac operator $D_M$ of the
boundary satisfies \cite[Theorem 6]{HM3}
\begin{eqnarray}\label{hmz}
\lambda_1 \geq \frac {n+1}{2} \inf_M H,
\end{eqnarray}
where $H$ denotes the (inward) mean curvature of $M$, assumed to be positive. Equality holds
if and only if $H$ is constant and every eigenspinor associated with $\lambda_1$ is the restriction
to $M$ of a parallel spinor on $N$ (and so $N$  is a Ricci-flat manifold). As an application of the limiting
case, they gave an elementary proof of the famous Alexandrov theorem: The only
compact embedded hypersurfaces of constant mean curvature in the Euclidean space are the round spheres.

\noindent The eigenvalue estimate (\ref{hmz}) has also led to several rigidity theorems \cite{HM1}. For example, if the boundary of a manifold is assumed to carry a Killing spinor (called extrinsic), then under some curvatures assumptions, the equality in \eqref{hmz} is attained and therefore the Killing spinor (which is the corresponding eigenspinor) comes from a parallel spinor on the whole manifold. As a consequence, the boundary has to be totally umbilical. This yields the following result: A complete Ricci-flat Riemannian manifold of dimension at least $3$ whose mean-convex boundary is isometric to the round sphere is a flat disc \cite[Cor. 6]{HM1}.

\noindent In a more general context, S.~Raulot assumed in \cite{Ra} the existence on the boundary of a solution of the Dirac equation (instead of an extrinsic Killing spinor), i.e. a spinor field $\varphi$ satisfying $D_M \varphi=\frac{n+1}{2}H_0\varphi$ for some particular function $H_0$. He showed that the boundary has to be connected and the spinor $\varphi$ is the restriction of a parallel spinor. In this case, the solution $\varphi$ satisfies an overdetermined differential equation involving the second fundamental form. As an application, he proved that if the whole manifold has a vanishing sectional curvature along the boundary (assumed to be simply connected), it has to be flat.

\noindent Recently, O.~Hijazi and S.~Montiel \cite{HM} covered some of the previous results by establishing an {\it integral inequality} relating the Dirac operator on the boundary to the mean curvature. Indeed they proved that, in general, if the mean curvature $H$ of the boundary is positive, then any spinor field $\varphi\in{\rm {\bf {S}}}:= \Sigma N\vert_M$ satisfies the inequality
\begin{equation}\label{hij-mon}
0\leq \int_M\Big(\frac{1}{H}|{\rm\bf D}_{\rm\bf{S}}\varphi|^2-\frac{(n+1)^2}{4}H|\varphi|^2\Big)dv,
\end{equation}
where $dv$ the volume element on $M$ and ${\rm\bf D}_{\rm\bf{S}}$
is the Dirac operator defined on ${\bf{S}}$: depending on the
dimension of the manifold $M$, the operator ${\rm\bf
D}_{\rm\bf{S}}$ coincides with $D_M$ or with the double copy
$D_M\oplus -D_M$, see Section \ref{sec:1} for more details.
Moreover, the equality holds in (\ref{hij-mon}) if and only if
there exists two parallel spinor fields $\psi,\theta\in
\Gamma(\Sigma N)$ such that $P_+\varphi=P_+\psi$ and
$P_-\varphi=P_-\theta$ on the boundary. Here the operators
$P_{\pm}$ are the orthogonal projections of $L^2(\textbf{S})$ onto
the $\pm 1$-eigenspaces corresponding to the $\pm 1$-eigenvalues
of the endomorphism $i\nu$ on $\textbf{S}$, where $\nu$ is the
inward unit vector field normal to the boundary.

\noindent The proof of Inequality (\ref{hij-mon}) relies mainly on the use of the spinorial Reilly formula applied to an appropriate spinor field which is a solution of a boundary value problem (see \cite{HM2} for a survey on those).
As an application of this inequality, O.~Hijazi and S.~Montiel proved that if $M$ admits an isometric and isospin immersion as an hypersurface into another spin Riemannian
manifold $\widetilde N$ carrying a parallel spinor field with
mean curvature $\widetilde H$, then the inequality
\begin{eqnarray}\label{STT}
\int_M H dv\leq \int_M\frac{\widetilde H^2}{H} dv,
\end{eqnarray}
holds \cite[Theorem 2]{HM}.
Moreover, the equality is realized if and only if the shape operators of both immersions are the same. In the special case where $\widetilde N = \mathbb R^{n+2}$, equality in (\ref{STT}) implies that $N$ is
a Euclidean domain and the immersion is congruent to the embedding of $M$ into
$N$ as its boundary. Inequality (\ref{STT}) implies also the Positive Mass Theorem in the spin case. Another application of interest is that the first eigenvalue of the Dirac operator of the boundary $M$ corresponding to the conformal metric with conformal factor $H^2$ is at least $\frac n2$ and equality holds if and only if there exists a non-trivial parallel spinor field on $N$ \cite[Theorem 1]{HM}.

\noindent In this paper, we are interested in studying the class of Riemannian spin manifolds whose
boundary is foliated by a unit vector field $\xi$ (we will call those foliations Riemannian
flows, see Section  \ref{sec:1}).
We aim at establishing rigidity results on such manifolds
when looking at solutions of the basic Dirac equation on the boundary, that is, a
spinor field $\varphi$ satisfying the equation
\begin{eqnarray}\label{D-E}
D_b\varphi=\frac{n+1}{2}H_0\varphi,
\end{eqnarray}
where $H_0$ is a basic function defined on the boundary and $D_b$ is the basic Dirac operator, see \cite{GK1,GK2} and references therein for the study of its transversal ellipticity. Here the word ``basic" refers to transverse objects that are constant along the leaves of the flow.

\noindent First, we start with studying solutions of the equation
(\ref{D-E}). In fact, we prove that  any solution $\varphi$
satisfies an integral inequality derived from (\ref{hij-mon}) (see
Theorem \ref{thm:5}) which involves the geometric data of the flow
and the boundary, such as the O'Neill tensor and the mean
curvature. By requiring the positivity of the $\xi$-direction
principal curvature, we show that the equality case of this
integral inequality restricts the geometry of the flow. In fact,
it has to be a local product (that is, the O'Neill tensor
vanishes) and the vector $\xi$ belongs to the kernel of the
Weingarten map (see Theorem \ref{equality-case}). In this case,
the spinor $\varphi$ (resp. $\varphi + \xi\cdot_M\varphi$) can be
extended to a parallel spinor on the whole manifold if $n$ is even
(resp. $n$ is odd). Here ``$\cdot_M$'' denotes the Clifford
multiplication on the boundary $M$.

\noindent In a second step and since basic Killing spinors (see e.g. \cite[eq. (1.4.5) p.37]{GH}) are particular solutions of the basic Dirac equation, we characterize the geometry of the flow carrying such spinors. By applying the equality case of the integral inequality, we prove the following rigidity result:

\begin{thm}\label{thm:6}
Let $(N,g)$ be an $(n+2)$-Riemannian spin manifold of non-negative scalar curvature with connected boundary $M$ of positive mean curvature $H$. Assume that $M$ is endowed with a minimal Riemannian flow carrying a maximal number of basic Killing spinors of constant $-\frac{1}{2}$  (resp. a maximal number of basic Killing spinors of constants $-\frac{1}{2}$ and $\frac{1}{2}$) if $n$ is even (resp. if $n$ is odd).
If the inequality $\frac{n}{n+1}+\frac{1}{n+1}[\frac{n}{2}]^\frac{1}{2}|\Omega|\leq H$ holds, then $N$ is isometric to the quotient $\lquot{\mathbb{R}\times B}{\Gamma}$ for some fixed-point-free cocompact discrete subgroup $\Gamma\subset\mathbb{R}\times\mathrm{SO}_{n+1}$, where $B$ is the closed unit ball in $\mathbb{R}^{n+1}$. 
In case $n$ is even, $\Gamma\cong\mathbb{Z}$ acts via $\left(k,(t,x)\right)\mapsto\left(t+ka,A^kx\right)$ for some $(a,A)\in\mathbb{R}^\times\times\mathrm{SO}_{n+1}$.
\end{thm}

\noindent Here $\Omega$ denotes the $2$-form associated to the O'Neill tensor field.
In the last part of the paper, we combine our result with an inequality by B.-Y. Chen \cite{Ch} in order to deduce a sufficient condition for having a non-trivial flow, i.e. a flow that is not a local product (see Theorem \ref{final}).

\section{Riemannian flows and manifolds with boundary} \label{sec:1}

\noindent In this section, we briefly recall some preliminaries on spin Riemannian flows (see \cite{T}, \cite{GH}). We then describe the geometric setting on manifolds with boundary. For more details, we refer to \cite{LM}, \cite{F}, \cite{G} and \cite{BHMMM}.
\subsection{Spin Riemannian flow}
\noindent Let $(M^{n+1},g)$ be a Riemannian manifold of dimension $n+1$ endowed with a Riemannian flow given by a unit vector field $\xi$. That is, the vector $\xi$ defines a one-dimensional foliation of $M$ given by its integral curves satisfying the bundle-like condition $\mathcal{L}_{\xi} g_{|_{\xi^\perp}}=0$ (see \cite{Rei}). Here $\mathcal{L}_\xi$ denotes the Lie derivative with respect to the vector field $\xi$. Equivalently, this means that the endomorphism $h$ defined by $h:= \nabla^M\xi:TM\rightarrow TM$ and known as the O'Neill tensor of the flow \cite{ON}, restricts to a skew-symmetric tensor field on the normal bundle $Q=\xi^\perp$. In particular, this gives rise to a unique metric connection $\nabla$ on the normal bundle, called transversal Levi-Civita connection, with respect to the induced metric. It is defined, for any section $Y \in \Gamma(Q)$, by
$$
\nabla _{X} Y =:
\left\{\begin{array}{ll}
\pi [X,Y]&\textrm {if $X=\xi$},\\\\
\pi (\nabla_{X}^{M}Y)&\textrm{if  $X\perp \xi$},
\end{array}\right.
$$
where $\pi:TM\rightarrow Q$ denotes the orthogonal projection \cite{T}. On the other hand, one can easily check that the corresponding Levi-Civita connections on $M$ and $Q$ are related for all sections $Z,W$ in $\Gamma(Q)$ via the Gauss-type formulas \cite{GH}:
\begin{equation*}
\left\{\begin{array}{ll}
\nabla^M_Z W=\nabla_Z W-g(h(Z),W)\xi, &\textrm {}\\\\
\nabla^M_\xi Z=\nabla_\xi Z+h(Z)-\kappa(Z)\xi,&\textrm {}
\end{array}\right.
\end{equation*}
where $\kappa:=\nabla^M_\xi \xi$ is the mean curvature of the flow.

\noindent From now on, we assume that $M$ is a spin manifold.
Since the tangent bundle of $M$ splits orthogonally as
$TM=\R\xi\oplus Q$, the pull-back of the spin structure on $M$
induces a spin structure on the normal bundle $Q$ (see \cite{B}).
In this case, the spinor bundle $\Sigma M$ is canonically
identified with the spinor bundle of $Q,$ denoted by $\Sigma Q,$
for $n$ even and with the direct sum $\Sigma Q\oplus \Sigma Q$ for
$n$ odd. The Clifford multiplications ``$\cdot_M$'' in $\Sigma M$
and ``$\cdot_Q$'' in $\Sigma Q$ are identified for any section
$Z\in \Gamma(Q)$ and spinor field $\varphi\in \Gamma(\Sigma Q)$ as
follows:
\begin{equation*}
\left\{\begin{array}{ll}
Z\cdot_M\varphi=Z \cdot_Q\varphi, &\textrm {for $n$ even}\\\\
Z\cdot_M\xi\cdot_M\varphi=(Z\cdot_Q\oplus-Z\cdot_Q)\varphi,&\textrm {for $n$ odd}.
\end{array}\right.
\end{equation*}

\noindent Therefore, by using of the above identification, one can relate the spinorial Levi-Civita connection $\nabla^M$ on $\Sigma M$ with the one on $\Sigma Q$ by \cite[eq. (2.4.7)]{GH}
\begin{equation}\label{eq:oneillspin}
\left\{\begin{array}{ll}\nabla_\xi^M\varphi&= \nabla_\xi\varphi+\frac{1}{2}\Omega\cdot_M\varphi+\frac{1}{2}\xi\cdot_M\kappa\cdot_M\varphi,\\
 &\\
\nabla_Z^M\varphi&= \nabla_Z\varphi+\frac{1}{2}\xi\cdot_M h(Z)\cdot_M\varphi,\end{array}\right.
\end{equation}
where $\Omega$ is the $2$-form associated to the tensor $h$ defined for all $Y,Z\in \Gamma(Q)$ by $\Omega(Y,Z)=g(h(Y),Z).$

\noindent Throughout this paper, we will consider Riemannian flows with basic mean curvature $\kappa$, i.e., the mean curvature of the flow satisfies $\nabla_\xi\kappa=0.$ The basic Dirac operator $D_b $ is a first order differential operator  (see \cite{GK1} and \cite{GK2})  defined on the set of basic sections (sections of the spinor bundle $\Sigma Q$ satisfying $\nabla_\xi\varphi=0$) by
$$D_b=\sum_{i=1}^n e_i\cdot_Q\nabla_{e_i}-\frac{1}{2}\kappa\cdot_Q,$$
where $\{e_i\}_{i=1,\cdots,n}$ is a local orthonormal frame of
$\Gamma(Q).$ Recall that the basic Dirac operator preserves the
set of basic sections and is a transversally elliptic and
essentially self-adjoint, if $M$ is compact. Therefore, by the
spectral theory of transversal elliptic operators, it is a
Fredholm operator and  has a discrete spectrum \cite{K1,K2}.

\noindent As a direct consequence from Equations \eqref{eq:oneillspin}, the transverse Levi-Civita connection commutes with the Clifford action of $\xi,$ that is $\nabla_X(\xi\cdot_M\varphi)=\xi\cdot_M\nabla_X\varphi$ for any spinor field $\varphi\in\Gamma(\Sigma Q)$ and $X\in \Gamma(TM).$ In particular, this means that the spinor field $\xi\cdot_M\varphi$ is basic if and only if $\varphi$ is basic. Thus, the basic Dirac operator shares a fundamental property with the Clifford action of $\xi$ on basic spinors. In fact, for $n$ even (resp. $n$ odd) and for any basic spinor field $\varphi$, we have
\begin{equation}\label{eq:commute}
D_b(\xi\cdot_M \varphi)=-\xi\cdot_M D_b \varphi\,\, \,\,({\rm resp.}\,\, D_b(\xi\cdot_M \varphi)=\xi\cdot_M D_b \varphi).
\end{equation}
Here we point out that, for $n$ odd, the identity in \eqref{eq:commute} is equivalent to saying that $(D_b\oplus -D_b)(\xi\cdot_M\varphi)=-\xi\cdot_M (D_b\oplus -D_b)\varphi,$ since for any $\varphi\in \Sigma Q\simeq \Sigma^+M,$ we have
$$-D_b(\xi\cdot_M\varphi)=(D_b\oplus-D_b)(\xi\cdot_M\varphi)=-\xi\cdot_M(D_b\oplus-D_b)\varphi=-\xi\cdot_M D_b\varphi.$$
Hence, the spectrum of the basic Dirac operator $D_b$ (resp.
$D_b\oplus -D_b$) is symmetric with respect to zero. We finish
this part by stating a fundamental relation between the Dirac
operator on $M$ and the basic Dirac operator. We have,
\begin{equation}\label{eq:relationdirac}
\left\{\begin{array}{ll}
D_M=D_b-\frac{1}{2}\xi\cdot_M\Omega\cdot_M, &\textrm {for $n$ even}\\\\
D_M=\xi\cdot_M(D_b\oplus -D_b)-\frac{1}{2}\xi\cdot_M\Omega\cdot_M, &\textrm {for $n$ odd}.
\end{array}\right.
\end{equation}

\subsection{Riemannian spin manifolds with boundary}
\noindent In the following, we will review some well-known facts
on spin manifolds with boundary, see e.g. \cite{HM1,HM2,HM3} and
references therein. Let $(N^{n+2},g)$ be a Riemannian spin
manifold of dimension $n+2$ with smooth boundary $M=\partial N.$
The existence of the (inward) unit vector field $\nu$ normal to
the boundary defines in a natural way a spin structure on $M$
induced from the one on $N$. This gives rise to the two spinor
bundles on the boundary, the intrinsic bundle $\Sigma M$ and the
extrinsic one  $\textbf{S}=\Sigma N|_{M}.$ Each bundle is  endowed
with the spinorial Levi-Civita connection, the Clifford
multiplication and the Dirac operator. The extrinsic data are
related to the ones on $N$ by the formulas:
\begin{eqnarray}\label{eq:gauss}
X\cdot_\textbf{S} \varphi &=&X\cdot \nu\cdot\varphi\nonumber\\
\nabla^N_X\varphi &=&\nabla^\textbf{S}_X \varphi +\frac{1}{2}A(X)\cdot_\textbf{S}\varphi\\
\textbf{D}_{\textbf{S}}\varphi&=&\frac{n+1}{2}H\varphi-\nu\cdot D_N\varphi-\nabla^N_\nu\varphi \nonumber,
\end{eqnarray}
where $``\cdot"$ is the Clifford multiplication on $N$, the tensor $A$ is the Weingarten map given for all $X\in \Gamma(TM)$ by $A(X)=-\nabla^N_X\nu,$ the spinor field $\varphi$ is a section in $\textbf{S}$ and $H=\frac{1}{n+1}{\rm Trace}(A)$ is the mean curvature of $M.$ The operator $\textbf{D}_\textbf{S},$ called the extrinsic Dirac operator, acts on sections on $\textbf{S}$ as $\textbf{D}_\textbf{S}=\sum_{i=1}^{n+1} e_i\cdot_\textbf{S}\nabla^{\textbf{S}}_{e_i},$ where $\{e_1,\cdots,e_{n+1}\}$ is a local orthonomal frame of $TM.$

\noindent On the other hand, the extrinsic spinor bundle can be identified with the intrinsic one in a canonical way depending on the dimension of $N.$ Namely, if $n$ is odd, the tuple $(\textbf{S},``\cdot_\textbf{S}",\nabla^\textbf{S},\textbf{D}_\textbf{S})$ can be identified to $(\Sigma M,``\cdot_{M}",\nabla^{M},D_{M})$ whereas for $n$ even it can be identified to $(\Sigma M\oplus \Sigma M,``\cdot_{M}\oplus -\cdot_{M}",\nabla^{M}\oplus \nabla^{M},D_{M}\oplus -D_{M}).$ Moreover, using the first two equations in \eqref{eq:gauss} and the Gauss formula, one can prove that the following relations hold for all $X,Y\in \Gamma(TM)$,
\begin{equation} \label{eq:compa}
\left\{\begin{array}{ll}
\nabla^\textbf{S}_X(Y\cdot)=\nabla^M_X Y\cdot+Y\cdot\nabla^\textbf{S}_X,\\\\
\nabla^\textbf{S}_X(\nu\cdot)=\nu\cdot\nabla_X^\textbf{S}
\end{array}\right.
\end{equation}
and that,
\begin{equation}\label{eq:espacepropre}
\textbf{D}_\textbf{S}(\nu\cdot)=-\nu\cdot\textbf{D}_\textbf{S}.
\end{equation}
Equality \eqref{eq:espacepropre} means that the spectrum of $\textbf{D}_\textbf{S}$ is symmetric with respect to zero and if $n$ is even the Dirac operator on $M$ commutes with the action of $\nu$, that is,
\begin{eqnarray} \label{eq:comu}
D_M(\nu\cdot\Phi)=\nu\cdot D_M\Phi
\end{eqnarray}
for any spinor field $\Phi\in \Gamma(\Sigma M).$

\noindent We define the operators $P_{\pm}$ as being the orthogonal projection of $L^2(\textbf{S})$ onto the $\pm 1$-eigenspaces corresponding to the $\pm 1$-eigenvalues of the operator $i\nu$ on $\textbf{S},$ i.e. $i\nu\cdot P_{\pm}=\pm P_{\pm}.$ They satisfy $$P_{\pm}(X \cdot)=X\cdot P_{\mp} \ \ \ \text{and} \ \ \ \ P_\pm(\nu\cdot)=\nu\cdot P_\pm,$$
for all $X\in \Gamma(TM)$. This implies that $\textbf{D}_\textbf{S}P_{\pm}=P_{\mp}\textbf{D}_\textbf{S}.$ Now, we state the following lemma about the spectrum of the basic Dirac operator when $N$ is an $(n+2)$-dimensional Riemannian spin manifold whose  boundary $M$ carries a Riemannian flow.
\begin{lemma} \label{neven}
Let $N$ be an $(n+2)$-dimensional Riemannian spin manifold with boundary $M$. Assume that $n$ is even and $M$ carries a Riemannian flow given by a unit vector field $\xi$. Then, for any basic spinor field $\varphi$, we have
$$D_b(\xi\cdot\varphi)=-\xi\cdot D_b\varphi.$$
\end{lemma}

\noindent{\bf Proof.} Using the identification of the spinor
bundles for $n$ even, we have
\begin{eqnarray}\label{iden}
\Sigma Q\oplus \Sigma Q\simeq \Sigma M\oplus\Sigma M\simeq
\textbf{S}.
\end{eqnarray}
Therefore, we can think of any $\varphi\in \Gamma(\Sigma Q)$ as a
section in one subbundle of $\textbf{S}$, say $\textbf{S}^+.$ This
means that $\xi\cdot\varphi$ is a section in $\textbf{S}^-$. In
order to prove the lemma, we need first to show that
$\xi\cdot\varphi$ is a basic spinor field. Indeed, using the first
equation in \eqref{eq:oneillspin} and the identification of the
Clifford multiplications in \eqref{eq:gauss}, we write
\begin{eqnarray*}
\nabla_\xi(\xi\cdot\varphi)& = &\nabla_\xi^M(\xi\cdot\varphi)-\frac{1}{2}\Omega\cdot(\xi\cdot\varphi)-\frac{1}{2}\xi\cdot\kappa\cdot (\xi\cdot\varphi)\\
&\bui{=}{\eqref{eq:gauss}}&\nabla^N_\xi(\xi\cdot\varphi)-\frac{1}{2}A\xi\cdot\nu\cdot\xi\cdot\varphi-\frac{1}{2}\xi\cdot\Omega\cdot\varphi-\frac{1}{2}\kappa\cdot\varphi\\
&=&\nabla^N_\xi\xi\cdot\varphi+\xi\cdot\nabla^N_\xi\varphi+\frac{1}{2}A\xi\cdot\xi\cdot\nu\cdot\varphi-\frac{1}{2}\xi\cdot\Omega\cdot\varphi-\frac{1}{2}\kappa\cdot\varphi\\
&\bui{=}{\eqref{eq:gauss}}&\nabla^M_\xi\xi\cdot\varphi+g(A\xi,\xi)\nu\cdot\varphi+\xi\cdot\nabla^M_\xi\varphi+\frac{1}{2}\xi\cdot A\xi\cdot\nu\cdot\varphi\\&&+\frac{1}{2}A\xi\cdot\xi\cdot\nu\cdot\varphi-\frac{1}{2}\xi\cdot\Omega\cdot\varphi-\frac{1}{2}\kappa\cdot\varphi\bui{=}{\eqref{eq:oneillspin}}\xi\cdot\nabla_\xi\varphi,
\end{eqnarray*}
which vanishes since $\varphi$ is a basic spinor field. Now using Equation \eqref{eq:relationdirac}, we have
\begin{eqnarray}\label{eq:dirac1}
D_M(\xi\cdot\varphi)&=&D_b(\xi\cdot\varphi)-\frac{1}{2}\xi\cdot_M\Omega\cdot_M(\xi\cdot\varphi)\nonumber\\
&=&D_b(\xi\cdot\varphi)+\frac{1}{2}\xi\cdot\nu\cdot\Omega\cdot\xi\cdot\varphi\nonumber\\
&=&D_b(\xi\cdot\varphi)+\frac{1}{2}\nu\cdot\Omega\cdot\varphi.
\end{eqnarray}

\noindent On the other hand and again by the identification between the Clifford multiplications, we compute
\begin{eqnarray} \label{eq:dirac2}
D_M(\xi\cdot\varphi)&=&D_M(\nu\cdot\xi\cdot\nu\cdot\varphi)
=D_M(\nu\cdot(\xi\cdot_M\varphi))\nonumber\\
&\bui{=}{\eqref{eq:comu}}& \nu\cdot D_M(\xi\cdot_M\varphi)\nonumber\\
&\bui{=}{\eqref{eq:relationdirac}}&\nu\cdot (D_b(\xi\cdot_M\varphi)-\frac{1}{2}\xi\cdot_M\Omega\cdot_M(\xi\cdot_M\varphi))\nonumber\\
&\bui{=}{\eqref{eq:commute}}&\nu\cdot (-\xi\cdot_M D_b\varphi+\frac{1}{2}\Omega\cdot_M\varphi)\nonumber\\
&=&-\xi\cdot D_b\varphi+\frac{1}{2}\nu\cdot\Omega\cdot\varphi.
\end{eqnarray}

\noindent Finally comparing Equations \eqref{eq:dirac1} and \eqref{eq:dirac2}, we deduce the desired identity.
\hfill$\square$\\

This lemma has the important consequence that if $\varphi$ is an eigenspinor for the basic Dirac operator associated to an eigenvalue $\lambda$, then the spinor $\xi\cdot\varphi$ is an eigenspinor associated to the eigenvalue $-\lambda$. 

\section{Foliated manifolds with boundary}\label{mainn}

\noindent In this section, we will deal with a Riemannian spin manifold whose boundary carries a Riemannian flow given by a unit vector field $\xi$. We will assume, after restricting the spin structure to the normal bundle, the existence of a spinor field $\varphi$ which is a solution for the basic Dirac equation. We will prove that the solution $\varphi$ satisfies an integral inequality involving the geometric data of the flow and the boundary.  When the equality case in this integral inequality is attained, we will show that the spinor solution $\varphi$ (resp. $\varphi + \xi\cdot_M \varphi$) is the restriction of a parallel spinor on the whole manifold $N$ if $n$ is even (resp. if $n$ is odd). The key point of our results is the integral inequality (\ref{hij-mon}) established by O. Hijazi and S. Montiel in \cite[Prop. 9]{HM}.

\noindent Next, we state the two main theorems of this section:
\begin{thm} \label{thm:5} Let $N$ be an $(n+2)$-dimensional  compact  Riemannian spin manifold with non-negative scalar curvature, whose boundary hypersurface $M$ has a positive mean curvature $H$ and is endowed with a Riemannian flow. Assume that there exists a spinor field $\varphi$ such that $D_b\varphi=\frac{n+1}{2}H_0\varphi,$ where $H_0$ is a positive basic function. Then, we have
\begin{eqnarray}\label{eq:estifoli}
0\leq \int_M\frac{1}{H}\big((n+1)^2 H_0^2|\varphi|^2 +
|\Omega\cdot_M\varphi|^2 -(n+1)^2 H^2|\varphi|^2 \big)dv.
\end{eqnarray}
\end{thm}

\noindent The equality case in Inequality \eqref{eq:estifoli} is now characterized by:

\begin{thm}\label{equality-case}
Under the same conditions as Theorem \ref{thm:5} and if we assume that $g(A(\xi),\xi)\geq 0$, then equality holds in (\ref{eq:estifoli}) if and only if $h=0$ (that is the flow is a local product) and $H_0=H.$ In this case, we get that $A(\xi)=0$ and the spinors $\varphi$ and $\xi\cdot\varphi$ are respectively the restrictions of parallel spinors on $N$ if $n$ is even, and if $n$ is odd the spinor $\varphi + \xi \cdot_M \varphi$ is the restriction of a parallel spinor on $N$.
\end{thm}

\noindent Since the identification between the data (the spinor bundles, the Clifford multiplications, the Dirac operators,...) of the manifold, the boundary and the normal bundle of the flow depend on the parity of $n$, we will distinguish two cases to prove Theorem \ref{thm:5} and Theorem \ref{equality-case}: the even-dimensional case and the odd-dimensional case.

\subsection{The even-dimensional case}

\noindent In this subsection, we will prove Theorem \ref{thm:5} and Theorem \ref{equality-case} for $n$ even.
We start by proving Theorem \ref{thm:5}.

\noindent{\bf Proof of Theorem \ref{thm:5}.} Since the manifold
$N$ is spin, the manifold $M$ and the normal bundle $Q$ of the
flow are also spin. Because $n$ is even, the spinor bundle of $Q$
is identified with the one of $M,$ which is also identified with
one subbundle of $\textbf{S}$, see \eqref{iden}. Therefore we can
think of any $\varphi\in \Gamma(\Sigma Q)$ as a section in one
subbundle of $\textbf{S}$, say $\textbf{S}^+.$ With the help of
Equation \eqref{eq:relationdirac}, we can say that
\begin{equation}\label{eq:1}
{\textbf
D}_{\textbf{S}}\varphi=D_M\varphi=\frac{n+1}{2}H_0\varphi-\frac{1}{2}\xi\cdot_M\Omega\cdot_M\varphi
\end{equation}
and
\begin{equation}\label{eq:1'}
{\textbf D}_{\textbf{S}} (\xi\cdot\varphi) = -D_M
(\xi\cdot\varphi) \bui{=}{\eqref{eq:dirac2}}  \frac{n+1}{2} H_0
\xi\cdot\varphi - \frac{1}{2} \nu \cdot \Omega\cdot\varphi .
\end{equation}
We then have
\begin{eqnarray*}
\vert {\textbf D}_{\textbf{S}} \varphi\vert^2 &=&
\frac{(n+1)^2}{4} H_0^2 \vert \varphi\vert^2 +\frac 14 \vert
\Omega\cdot_M\varphi\vert^2 - \frac{n+1}{2} H_0 \langle \varphi,
\xi\cdot_M \Omega\cdot_M\varphi \rangle
\end{eqnarray*}
and
\begin{eqnarray*}
\vert {\textbf D}_{\textbf{S}}(\xi\cdot\varphi)\vert^2 &=&
\frac{(n+1)^2}{4} H_0^2 \vert \varphi\vert^2 +\frac 14 \vert
\Omega\cdot_M\varphi\vert^2  + \frac{n+1}{2} H_0 \langle\varphi,
\xi\cdot_M \Omega\cdot_M\varphi \rangle .
\end{eqnarray*}
Thus, by applying Inequality (\ref{hij-mon}) to $\varphi$ and
$\xi.\varphi$ and summing the obtained inequalities, we find the
desired result. \hfill$\square$

\noindent Before proving Theorem \ref{equality-case}, we start by establishing a crucial lemma about
the characterization of the equality case of \eqref{eq:estifoli}.
\begin{lemma}
Under the same conditions as Theorem \ref{thm:5} and if equality
holds in (\ref{eq:estifoli}), we have
\begin{eqnarray} \label{eq:4}
h(X)\cdot_M \varphi+g(A(\xi),X)\frac{H_0}{H}\varphi-\frac{1}{(n+1)H}g(A(\xi),X)\xi\cdot_M\Omega\cdot_M \varphi\nonumber\\
=-\frac{1}{(n+1)H}A(X)\cdot_M\Omega\cdot_M \varphi,
\end{eqnarray}
for all $X\in \Gamma(TM).$
\end{lemma}

\noindent{\bf Proof.}  Assume that equality holds in
(\ref{eq:estifoli}). We know that the optimality is characterized
by the existence of two parallel spinors $\psi,\theta$ on $N$ such
that $P_+\varphi=P_+\psi$ and $P_-\varphi=P_-\theta$ on $M$. By
taking ${\textbf D}_{\textbf{S}}$ on both sides and using the last
equation in \eqref{eq:gauss} for the spinors $\psi$ and $\theta$,
we get $P_-({\textbf D}_{\textbf{S}}\varphi)=\frac{n+1}{2}H
P_-\psi$ and $P_+({\textbf D}_{\textbf{S}}\varphi)=\frac{n+1}{2}H
P_+\theta.$

\noindent There exists two parallel spinors $\Psi,\Theta$ such
that the following holds on $M$:
$$P_+(\xi\cdot\varphi)=P_+\Psi\ \ \ \ \text{and}\ \ \ \  P_-(\xi\cdot\varphi)=P_-\Theta.$$
As before, this also means that
$$P_-({\textbf D}_{\textbf{S}}(\xi\cdot\varphi))=\frac{n+1}{2}H P_-\Psi \ \ \ \text{and}\ \ \  P_+({\textbf D}_{\textbf{S}}(\xi\cdot\varphi))=\frac{n+1}{2}H P_+\Theta.$$

\noindent Now, differentiating the equation $\xi\cdot P_-\theta=P_+\Psi$ along any vector field $X$ in $\Gamma(TM)$ and using the first formula in \eqref{eq:compa} with the fact that  $\nabla^{\textbf{S}}$ commutes with $P_\pm$, we get
$$h(X)\cdot P_-\theta+\xi\cdot P_-(\nabla^{\textbf{S}}_X\theta)=P_+(\nabla^{\textbf{S}}_X\Psi).$$
Since the spinors $\theta$ and $\Psi$ are parallel, we deduce from \eqref{eq:gauss} that the above equation can be reduced to
\begin{equation}\label{eq:step1}
h(X)\cdot P_-\theta-\frac{1}{2}\xi\cdot A(X)\cdot \nu\cdot P_+\theta=-\frac{1}{2}A(X)\cdot\nu\cdot P_-\Psi.
\end{equation}

\noindent In the same way and by differentiating the equation $\xi\cdot P_+\psi=P_-\Theta,$  we can derive the relation
$$h(X)\cdot P_+\psi-\frac{1}{2}\xi\cdot A(X)\cdot \nu\cdot P_-\psi=-\frac{1}{2}A(X)\cdot\nu\cdot P_+\Theta.$$
\noindent On the other hand
\begin{eqnarray*}
{\textbf
D}_{\textbf{S}}(\xi\cdot\varphi)\bui{=}{\eqref{eq:1'}}\frac{n+1}{2}
H_0 \xi\cdot\varphi - \frac{1}{2} \nu \cdot \Omega\cdot\varphi
\bui{=}{\eqref{eq:1}} (n+1)H_0\xi\cdot\varphi-\xi\cdot {\textbf
D}_{\textbf{S}}\varphi.
\end{eqnarray*}

\noindent Thus, by applying $P_\mp$ on both sides of the last equality we find that
\begin{equation}\label{eq:2}
HP_-\Psi=2H_0\xi\cdot P_+\varphi-H\xi\cdot P_+\theta,
\end{equation}
and,
$$HP_+\Theta=2H_0\xi\cdot P_-\varphi-H\xi\cdot P_-\psi.$$
In  Equation \eqref{eq:step1}, we replace $P_-\Psi$ by its value from \eqref{eq:2} and $P_+\theta$ by
\begin{eqnarray} \label{eq:18}
P_+\theta=\frac{2}{(n+1)H}P_+({\textbf D}_{\textbf{S}}\varphi)
\bui{=}{\eqref{eq:1}}\frac{H_0}{H}P_+\varphi-\frac{1}{(n+1)H}\xi\cdot\nu\cdot \Omega\cdot P_-\varphi,
\end{eqnarray}
to get  the following identity
\begin{eqnarray} \label{I1}
h(X)\cdot P_-\varphi &+& g(A(\xi),X)\frac{H_0}{H}\nu\cdot P_+\varphi -\frac{1}{(n+1)H}g(A(\xi),X)\xi\cdot\Omega\cdot P_-\varphi  \nonumber \\
&=&-\frac{1}{(n+1)H}A(X)\cdot\Omega\cdot P_-\varphi.
\end{eqnarray}
Also, we can prove that a similar identity as in (\ref{I1}) holds with the opposite sign. This gives after summing them Equation (\ref{eq:4}).
\hfill$\square$

\noindent We are now able to prove Theorem \ref{equality-case} for $n$ even.

\noindent{\bf Proof of Theorem \ref{equality-case}.} Assume that Equality holds in (\ref{eq:estifoli}) and let us  prove that $h=0$ and $A(\xi)=0$. First of all, multiplying \eqref{eq:4} by a vector $X$ and contracting over an orthonormal frame $\{\xi,e_1,\cdots,e_n\}$ of $TM$, we get that
\begin{equation}\label{eq:6}
\xi\cdot_M\kappa\cdot_M\varphi+\Omega\cdot_M\varphi+\frac{H_0}{H}A(\xi)\cdot_M\varphi-\frac{1}{(n+1)H}A(\xi)\cdot_M\xi\cdot_M\Omega\cdot_M \varphi=0,
\end{equation}
where we used the local expression of  $\Omega$ given by $\Omega=\frac{1}{2}\sum_{i=1}^n e_i\cdot_M h(e_i).$ Multiplying again Equation \eqref{eq:6} by $A(\xi)\cdot_M \xi\cdot_M$ and plugging the last term in \eqref{eq:6} by its value into the new equation, we find after a straightforward computation that
\begin{eqnarray}\label{eq:8}
-A(\xi)\cdot_M\kappa\cdot_M\varphi+\mathcal B \ \Omega\cdot_M\varphi+\frac{H_0}{H}|A(\xi)|^2\xi\cdot_M\varphi+(n+1)H_0A(\xi)\cdot_M\varphi\nonumber\\
+\Big((n+1)H+2g(A(\xi),\xi)\Big)\xi\cdot_M\kappa\cdot_M\varphi=0,
\end{eqnarray}
where we denote by $\mathcal B$ the term
$$\mathcal B=(n+1)H+2g(A(\xi),\xi)+\frac{1}{(n+1)H}|A(\xi)|^2.$$
\noindent We now claim that $\mathcal B\neq 0.$ In fact, we can
write  $\mathcal B =  \frac{1}{(n+1)H} \vert A(\xi) + (n+1)H
\xi\vert^2.$ If $\mathcal B = 0$, then $A(\xi) = -(n+1)H \xi$.
Plugging this in (\ref{eq:6}), we get $\xi\cdot_M\kappa\cdot_M
\varphi = H_0 (n+1)\xi\cdot_M\varphi$ and so  $\kappa\cdot_M
\varphi = H_0 (n+1)\varphi$. Taking the real part of the scalar
product of the last identity with $\varphi$, we get
$$
(n+1) H_0 \vert \varphi\vert^2 = 0.
$$
Since by assumption $H_0>0$, this gives $\varphi=0$. But the
zero-set of any basic solution $\varphi$ of \eqref{D-E} has dense
complement subset in $M$ by \cite[Main Theorem]{Baer97} since the
basic Dirac operator on basic sections only differs from the Dirac
operator $D_M$ by a zero-order-term, see \eqref{eq:relationdirac}.
Therefore, we obtain a contradiction and hence $\mathcal B \neq 0$
must hold. Now, we claim the following:
\begin{lemma} \label{calc}We have that
 $$H\mathcal I \mathcal L \cdot_M\kappa\cdot_M\varphi+H_0 \mathcal J\mathcal  L \cdot_M\varphi=0,$$
where
$$\mathcal I:=(n+1)H+g(A(\xi),\xi),$$
$$\mathcal J:=(n+1)Hg(A(\xi),\xi)+|A(\xi)|^2,$$
$$\mathcal K:=(n+1)H+2g(A(\xi),\xi),$$
$$\mathcal L:=A(\xi)-\mathcal K\xi.$$
\end{lemma}
Assume the lemma holds for the moment. First of all, we point out
that $\mathcal L$ is a non-vanishing vector field. Indeed, assume
that $A(\xi)=\mathcal K \xi$ then, $g(A(\xi),\xi)=-(n+1)H \xi$.
This implies that $A(\xi)=-(n+1)H\xi$ which leads to a
contradiction because $\mathcal B \ne 0$. Therefore, we deduce
that $H \mathcal I \kappa\cdot_M\varphi+H_0 \mathcal J \varphi=0$
and, and, since the zero-set of $\varphi$ has dense complement
subset in $M$, we obtain $\mathcal J=0$ and $\mathcal I\kappa=0.$
Notice that $\mathcal I$ cannot be zero because $\mathcal I =-
\frac {1}{(n+1)H} \mathcal J + \mathcal B = \mathcal B$. Thus
$\kappa=0$ and $A(\xi)=0$. Plugging this in  (\ref{eq:8}) gives
that $\Omega\cdot_M \varphi =0$. Finally, Equation \eqref{eq:4}
implies that $h=0$. \noindent In order to prove the last statement
of Theorem \ref{equality-case}, we proceed as in \cite[Thm.
2]{HM}. Using Equation \eqref{eq:18}, we have that
$$H_0P_+\varphi=H P_+\theta \quad\text{and}\quad H_0 P_-\varphi=HP_-\psi.$$
Applying ${\textbf D}_{\textbf{S}}$ on both sides of the two equalities yields after taking the sum that
$$dH_0\cdot_{\textbf{S}}\varphi+\frac{n+1}{2}H_0^2\varphi=\frac{H_0}{H}dH\cdot_{\textbf{S}}\varphi+\frac{n+1}{2}H^2\varphi.$$
Thus, the Hermitian product by $\varphi$ gives $H_0=H$ which means
that $P_+\varphi=P_+\theta$ and $P_-\varphi=P_-\psi$. Thus, we get
that $\varphi=\psi=\theta$ on $M$. The same can be done for
$\xi\cdot\varphi.$
\hfill$\square$\\\\
We still have to prove Lemma \ref{calc}. \\\\
{\bf Proof of Lemma \ref{calc}.} Since $\mathcal B \neq 0$, Equation (\ref{eq:8}) implies that
\begin{eqnarray}\label{omegaB}
\Omega\cdot_M \varphi = &\frac{1}{\mathcal B}& \Big (
A(\xi)\cdot_M \kappa\cdot_M \varphi - [(n+1)H + 2 g(A(\xi), \xi)]
\xi\cdot_M\kappa\cdot_M \varphi\nonumber \\ && - (n+1)H_0 A(\xi)
\cdot_M \varphi - \frac{H_0}{H} \vert A (\xi)\vert^2 \xi \cdot_M
\varphi \Big ).
\end{eqnarray}
Taking the Clifford multiplication of Equation (\ref{eq:4}) with $\xi$ and for $X = \xi$, we get
\begin{eqnarray}\label{eq:44}
\xi\cdot_M \kappa\cdot_M\varphi &+& g(A(\xi), \xi) \frac{H_0}{H}
\xi\cdot_M \varphi  \nonumber \\ && = \frac{1}{(n+1)H}
\Big[A(\xi)\cdot_M \xi \cdot_M+  g(A(\xi), \xi)\Big]\Omega\cdot_M
\varphi .
\end{eqnarray}
Combining Equations (\ref{eq:44}) and (\ref{omegaB}), we get
\begin{eqnarray*}
(n+1) H &\mathcal B& \Big [ \xi\cdot_M \kappa\cdot_M\varphi+
g(A(\xi), \xi) \frac{H_0}{H} \xi\cdot_M \varphi\Big]  \nonumber \\
&& =  \vert A(\xi)\vert^2 \xi\cdot_M \kappa\cdot_M \varphi -2
g(A(\xi), \xi) A(\xi)\cdot_M \kappa\cdot_M \varphi \\ && + (n+1)H
A(\xi)\cdot_M \kappa\cdot_M \varphi + 2 g(A(\xi), \xi)
A(\xi)\cdot_M \kappa\cdot_M \varphi \\ && - (n+1)H_0 \vert
A(\xi)\vert^2 \xi\cdot_M \varphi + 2 (n+1)H_0 g(A(\xi), \xi)
A(\xi) \cdot_M \varphi \\ && + \frac{H_0}{H} \vert A(\xi)\vert^2
A(\xi) \cdot_M \varphi + g(A(\xi), \xi) A(\xi) \cdot_M
\kappa\cdot_M \varphi \\ && - (n+1)Hg(A(\xi), \xi) \xi \cdot_M
\kappa\cdot_M \varphi - 2 g(A(\xi), \xi)^2 \xi \cdot_M
\kappa\cdot_M \varphi \\ && - (n+1)H_0  g(A(\xi), \xi) A(\xi)
\cdot_M \varphi - \frac{H_0}{H} \vert A (\xi)\vert^2 g(A(\xi), \xi) \xi\cdot_M\varphi.
\end{eqnarray*}
The last expression can be written as $C_1 \xi \cdot_M
\kappa\cdot_M \varphi + C_2  \xi\cdot_M \varphi = C_3 A(\xi)
\cdot_M \kappa\cdot_M \varphi + C_4 A(\xi)\cdot_M \varphi$. Let us
determine $C_1, C_2, C_3$ and $C_4$. For example, the coefficient
$C_1$ of  $\xi \cdot_M \kappa\cdot_M \varphi$ is given by
\begin{eqnarray*}
C_1 &=& (n+1) H \mathcal B - \vert A(\xi)\vert^2 + (n+1) H
g(A(\xi), \xi) +2 g(A(\xi), \xi)^2 \\&=& (n+1)^2H^2 +  3 (n+1) H
g(A(\xi), \xi)  +2 g(A(\xi), \xi)^2 \\ &=& \Big((n+1)H+
g(A(\xi),\xi)\Big)\Big((n+1)H+2g(A(\xi),\xi)\Big) \\
&=& \mathcal I \mathcal K .
\end{eqnarray*}
In a similar way, we can find
$$C_2 =\frac{H_0}{H}\Big((n+1)Hg(A(\xi),\xi)+|A(\xi)|^2\Big)\Big((n+1)H+2g(A(\xi),\xi)\Big)=\frac{H_0}{H} \mathcal J \mathcal K,$$
$$C_3 = (n+1)H+g(A(\xi),\xi)= \mathcal I,$$
$$C_4 = \frac{H_0}{H}\Big((n+1)Hg(A(\xi),\xi)+|A(\xi)|^2\Big)=\frac{H_0}{H}\mathcal J .$$
Finally, we have
$$
H \mathcal I \mathcal K \xi\cdot_M\kappa\cdot_M\varphi +
H_0\mathcal J \mathcal K\xi\cdot_M\varphi = H\mathcal I
A(\xi)\cdot_M \kappa\cdot_M\varphi +H_0 \mathcal J
A(\xi)\cdot_M\varphi,
$$
which is $H\mathcal I \mathcal L \cdot_M\kappa\cdot_M\varphi+H_0
\mathcal J\mathcal  L \cdot_M\varphi=0$. \hfill$\square$

\begin{remark}
We notice that when equality case in (\ref{eq:estifoli}) is
realized, the second fundamental form has two important
properties:
\begin{enumerate}
\item First, it commutes with the O'Neill tensor. Indeed, the
Hermitian product of \eqref{eq:4} with $\xi\cdot_M\varphi$
provides the vanishing of the real part of $\langle
A(X)\cdot_M\Omega\cdot_M\varphi,\xi\cdot_M\varphi\rangle$ for any
$X\in \Gamma(TM).$ Therefore, by taking again the Hermitian
product of \eqref{eq:4} with $A(X)\cdot_M\varphi$, we obtain
$g(h(X),A(X))=0$ which means that $h$ and $A$ commute. In
particular, the eigenvalues of $A$ are of multiplicity at least
two if $h$ is not zero. \item Second, it has two constant
principal curvatures $0$ and $-(n+1)H$. In fact, by taking the
Hermitian product of \eqref{eq:4} with $\varphi$ and plugging
Equation \eqref{omegaB}, we get after a straightforward
computation that $A^2(\xi)=-(n+1)HA(\xi).$
\end{enumerate}
\end{remark}

\subsection{The odd-dimensional case}

\noindent In this subsection, we will prove Theorem \ref{thm:5}
and Theorem \ref{equality-case} when $n$ is odd. The difference
between the proofs in both cases comes certainly from the
different identification of the spinor bundles, the Clifford
multiplications and the Dirac operators.

\noindent We recall from Section \ref{sec:1} the identifications
of the spinor bundles of $N,$ the boundary $M$ and the normal
bundle of the flow as
$$\Sigma Q\oplus \Sigma Q\simeq \Sigma M\simeq \textbf{S},$$
where in the first isomorphism, the Clifford multiplications are being identified as
$$(Z\cdot_Q\oplus -Z\cdot_Q)\Upsilon=Z\cdot_M\xi\cdot_M\Upsilon,$$
for any $Z\in\Gamma(Q),$ while in the second isomorphism
$X\cdot_M\Upsilon=X\cdot\nu\cdot\Upsilon$ for any spinor field
$\Upsilon\in \textbf{S}=\Sigma N|_M.$ From the fact that $n+1$ is
even, the action of $i\nu$ on $\textbf{S}$ is determined by the
action of the complex volume form $\omega$ of $\Sigma M$, that is
for any spinor $\Upsilon$ on $\textbf{S},$ we have
$i\nu\cdot\Upsilon=\omega\cdot\Upsilon=\bar{\Upsilon}$, where
$\bar{\Upsilon}=\Upsilon_+ -\Upsilon_{-}$ with $\Upsilon_\pm$ are
eigensections of $\omega$ corresponding to the eigenvalues $\pm
1.$ Thus, from the definition of the projections $P_\pm$, we
deduce that $P_\pm\Upsilon=\Upsilon_\pm$.

\noindent{\bf Proof of Theorem \ref{thm:5}.} As for the even case,
the proof is mainly based on the use of Inequality (\ref{hij-mon})
applied to the spinor field $\Upsilon=\varphi+\xi\cdot_M\varphi.$
Here $\varphi$ is considered as a section in $\Sigma Q\simeq
\Sigma^+M,$ i.e. $P_+\varphi=\varphi.$ Therefore, we compute
\begin{eqnarray} \label{eq:relations}
{\textbf D}_{\textbf{S}}\Upsilon &=&D_M\Upsilon =D_M\varphi+D_M(\xi\cdot_M\varphi)\nonumber\\
&\bui{=}{\eqref{eq:relationdirac}}& \frac{n+1}{2}H_0\xi\cdot_M\varphi-\frac{1}{2}\xi\cdot_M\Omega\cdot_M\varphi-\xi\cdot_M D_b(\xi\cdot_M\varphi)+\frac{1}{2}\Omega\cdot_M\varphi\nonumber\\
&\bui{=}{\eqref{eq:commute}}& \frac{n+1}{2}H_0\xi\cdot_M\varphi-\frac{1}{2}\xi\cdot_M\Omega\cdot_M\varphi+\frac{n+1}{2}H_0\varphi+\frac{1}{2}\Omega\cdot_M\varphi.
\end{eqnarray}
It is easy to check that $\langle \xi\cdot_M\Omega\cdot_M\varphi,
\varphi \rangle=0$, we then deduce that the norm of the spinor
field $|{\textbf D}_{\textbf{S}}\Upsilon|^2$ is equal to
$$|{\textbf D}_{\textbf{S}}\Upsilon|^2=\frac{(n+1)^2}{2}H_0^2|\varphi|^2+\frac{1}{2}|\Omega\cdot_M\varphi|^2.$$
Plugging the last equality into Inequality (\ref{hij-mon}) and using the fact that $|\Upsilon|^2=2|\varphi|^2$, we get Inequality (\ref{eq:estifoli}). \hfill$\square$

\noindent Before proving Theorem \ref{equality-case} and as we did for the even case, we will start
by stating and proving a lemma which characterizes the equality case in (\ref{eq:estifoli}).

\begin{lemma}
Under the same conditions as Theorem \ref{thm:5} and if equality holds in (\ref{eq:estifoli}) the spinor field $\Upsilon : =\varphi + \xi\cdot_M\varphi$ satisfies
\begin{eqnarray} \label{eq:4*}
h(X)\cdot_M \Upsilon+g(A(\xi),X)\frac{H_0}{H}\Upsilon-\frac{1}{(n+1)H}g(A(\xi),X)\xi\cdot_M\Omega\cdot_M \Upsilon\nonumber\\
=-\frac{1}{(n+1)H}A(X)\cdot_M\Omega\cdot_M \Upsilon,
\end{eqnarray}
for all $X\in \Gamma(TM).$
\end{lemma}

\noindent{\bf Proof.} Assume that we have the equality case in (\ref{eq:estifoli}). As before, the optimality is characterized by the existence of two parallel spinors $\Psi$ and $\Theta$ such that on $M$ we have
$$P_+\Upsilon=\varphi=P_+\Psi\quad\text{and}\quad P_-\Upsilon=\xi\cdot_M\varphi=P_-\Theta.$$
\noindent Differentiating the second equation with respect to any vector field $X\in \Gamma(TM),$ we find after using Equation \eqref{eq:compa} that
\begin{equation}\label{eq:compat}
h(X)\cdot\nu\cdot\varphi+\frac{i}{2}\xi\cdot A(X)\cdot\nu\cdot P_-\Psi=-\frac{1}{2}A(X)\cdot\nu\cdot P_+\Theta.
\end{equation}
On the other hand, we get from Equation \eqref{eq:relations} that
$${\textbf D}_{\textbf{S}}\varphi=\frac{n+1}{2}H_0\xi\cdot_M\varphi-\frac{1}{2}\xi\cdot_M\Omega\cdot_M\varphi=\frac{n+1}{2}HP_-\Psi,$$
and,
$${\textbf D}_{\textbf{S}}(\xi\cdot_M\varphi)=\frac{n+1}{2}H_0\varphi+\frac{1}{2}\Omega\cdot_M\varphi=\frac{n+1}{2}HP_+\Theta.$$
Hence, replacing the spinor fields $P_-\Psi$ and $P_+\Theta$ by their values from the above two equations into \eqref{eq:compat}, we find after a straightforward computation that
\begin{eqnarray} \label{part1}
h(X)\cdot_M\varphi+\frac{H_0}{H}g(A(X),\xi)\xi\cdot_M\varphi+\frac{1}{(n+1)H}A(X)\cdot_M\Omega\cdot_M\varphi=\nonumber\\
\frac{1}{(n+1)H}g(A(X),\xi)\xi\cdot_M\Omega\cdot_M\varphi.
\end{eqnarray}
Here we used the identification between the Clifford multiplications on $\Sigma M$ and $\Sigma N$ and the fact that $i\nu\cdot\varphi=\varphi$.
Now one can easily check that Equation (\ref{eq:4*}) holds for the spinor field $\Upsilon$
after taking the Clifford multiplication of (\ref{part1}) by the vector $\xi$.
\hfill$\square$

\noindent Now, we prove Theorem \ref{equality-case} for $n$ odd.

\noindent{\bf Proof of Theorem \ref{equality-case}.} We remark
that Equation (\ref{eq:4*}) is the same as Equation \eqref{eq:4}
when $n$ was even but the spinor $\varphi$ in (\ref{eq:4}) is
replaced by the spinor $\Upsilon= \varphi + \xi \cdot_M \varphi$
in (\ref{eq:4*}). So, we repeat the same technique applied in the
proof of Theorem \ref{equality-case} when $n$ was even.
\hfill$\square$

\begin{remark}
The condition taken in the equality case of Theorem \ref{thm:5}
can be refined by requiring that the sectional curvature on $N$
vanishes along planes in $TM$ containing $\xi$. Indeed, using the
Gauss-Codazzi equations the Ricci curvature on $M$ of the vector
$\xi$ is equal to
$${\rm Ric}_M\xi=\sum_{i=1}^n R^N(\xi,e_i)e_i-A^2(\xi)+(n+1)HA(\xi),$$
where $\{e_i\}$ is a local orthonormal frame of $\Gamma(Q)$. Since $g({\rm Ric}_M\xi,\xi)=|h|^2$ \cite{ON}, we deduce that $g(A(\xi),\xi)\geq 0.$
\end{remark}

\section{Rigidity results for the geometry of the manifold, its boundary and the Riemannian flow}

\noindent In this section, we will state various rigidity results
on manifolds whose boundary carries a particular solution of the
basic Dirac equation and the mean curvature is assumed to satisfy
some condition depending on the norm of the O'Neill tensor of the
flow.  These results can be seen as the foliated counterpart of
the ones in \cite{Ra}.

\noindent We will estimate the two terms in the integral
inequality (\ref{eq:estifoli}) involving the $2$-form $\Omega$ in
terms of its norm. Indeed, we have:

\begin{prop} \label{thm:rig2} Under the assumptions of Theorem \ref{thm:5}, the inequality
\begin{equation} \label{est}
0\leq\int_M\frac{1}{H}\left((n+1)^2H_0^2
+[\frac{n}{2}]|\Omega|^2-(n+1)^2 H^2\right)|\varphi|^2dv,
\end{equation} holds. Moreover, under the condition $(n+1)H_0
+[\frac{n}{2}]^{\frac{1}{2}}|\Omega| \le (n+1) H$, equality is
attained if and only if $h=0$ and $H_0=H$. In this case, we have
that $A(\xi)$=0 and the spinors $\varphi$ and $\xi\cdot\varphi$
are respectively the restrictions of parallel spinors on $N$ if
$n$ is even, and if $n$ is odd the spinor $\varphi + \xi \cdot_M
\varphi$ is the restriction of a parallel spinor on $N$.
\end{prop}

\noindent{\bf Proof.} Since the operator $i\Omega$ is Hermitian,
all its eigenvalues are real. Therefore, one can always find an
orthonormal frame $\{e_i\}$ of $TM$ such that
$$\Omega=\sum_{j=1}^{[\frac{n}{2}]}\lambda_j e_{2j-1}\wedge e_{2j}.$$
Here $\lambda_j$ are the eigenvalues of the operator $i\Omega$.
Thus, from the Cauchy-Schwarz inequality, we get
\begin{equation*}
|\Omega\cdot_M\varphi|\leq \sum_{j=1}^{[\frac{n}{2}]} |\lambda_j|
|\varphi|\leq [\frac{n}{2}]^\frac{1}{2}(\sum_{j=1}^{[\frac{n}{2}]}
|\lambda_j|^2)^\frac{1}{2}|\varphi|\leq
[\frac{n}{2}]^\frac{1}{2}|\Omega||\varphi|,
\end{equation*}
and we conclude \eqref{est} by using Theorem \ref{thm:5}. \\
\noindent The equality holds if and only if $(n+1)^2H_0^2
+[\frac{n}{2}]|\Omega|^2=(n+1)^2 H^2$. In this case, the condition
$(n+1)H_0 +[\frac{n}{2}]^{\frac{1}{2}}|\Omega| \le (n+1) H$ allows
us to deduce that $\Omega=0$ which yields from \eqref{eq:4} to
$h=0$ and $A(\xi)=0.$ When $n$ is odd, the same technique carries
again over by using \eqref{eq:4*}. \hfill$\square$
\begin{remark}
 Under the assumptions of Theorem \ref{thm:5} and the conditions $(n+1)^2H_0^2
+[\frac{n}{2}]|\Omega|^2\le (n+1)^2 H^2$ and $g(A(\xi),\xi) \ge
0$, equality holds in \eqref{eq:estifoli}. Then, by Theorem
\ref{equality-case}, we have $h=0$ and $A(\xi)=0$.
\end{remark}

\begin{remark}
It is an easy fact to see that if there exists a basic harmonic
spinor on the flow, that is $D_b\varphi=0$, the estimate
$$\int_M\frac{1}{H}\left([\frac{n}{2}]|\Omega|^2-\frac{(n+1)^2}{4}H^2\right)|\varphi|^2dv>0,$$
holds. However, the limiting case cannot be achieved since the
mean curvature is assumed to be positive.
\end{remark}

\noindent Next, we deduce an analogue result of the main theorem in \cite[Thm. 1]{Ra}:

\begin{cor} \label{thm:rig} Let $N$ be a compact spin Riemannian $(n+2)$-dimensional manifold with non-negative scalar curvature,
whose boundary hypersurface $M$ has positive mean curvature $H$
and is endowed with a Riemannian flow. Assume that there exist a
spinor field $\varphi$ such that
$D_b\varphi=\frac{n+1}{2}H_0\varphi,$ where $H_0$ is a positive
basic function with
$H_0+\frac{1}{n+1}[\frac{n}{2}]^\frac{1}{2}|\Omega|\leq H$. Then
the vector field $\xi$ is parallel on $M$ and $A(\xi)=0$.
Moreover, the spinors $\varphi$ and $\xi\cdot\varphi$ are
respectively the restrictions of parallel spinors on $N$ if $n$ is
even and if $n$ is odd, the spinor $\varphi + \xi \cdot_M \varphi$
is the restriction of a parallel spinor on $N$.
\end{cor}

\noindent In the following subsection, we will apply Corollary
\ref{thm:rig} for particular solutions of the basic Dirac
equation, namely the basic Killing spinors \cite{GH} in order to
prove Theorem \ref{thm:6}.

\subsection{Basic Killing spinors as solutions of the basic Dirac equation}

\noindent It is a standard fact that on a spin manifold, Killing
spinors are particular solutions of the Dirac equation. Those
spinors appeared in the limiting case of an eigenvalue estimate of
the Dirac operator (see e.g. \cite{F} and references therein) and
they are completely classified in \cite{B2}. In the following, we
will consider a Riemannian flow defined on the boundary of a spin
manifold and will assume that it carries a basic Killing spinor,
that is, spinor $\varphi$ satisfying $\nabla_X\varphi=\alpha
X\cdot_Q\varphi$ for all $X\in \Gamma(Q)$ and
$\nabla_\xi\varphi=0$ for a real number $\alpha$ (assumed to be
$\pm\frac{1}{2}$). We will prove Theorem \ref{thm:6} by
considering again the two cases: $n$ is even and $n$ is odd.

\noindent{\bf Proof of Theorem \ref{thm:6} for $n$ even.} Note
first that the existence of a basic Killing spinor $\varphi$ of
Killing constant $-\frac{1}{2}$ gives rise to a solution of the
basic Dirac equation  with $H_0=\frac{n}{n+1},$ since the flow is
minimal. Hence, under the assumption on the mean curvature, we
deduce from Corollary \ref{thm:rig} that there exists a maximal
number of parallel spinors on $N$ (which restrict to $\varphi$ and
$\xi\cdot\varphi$) and that $\nabla^M\xi=0$ and $A(\xi)=0$.
Moreover, we find from Equations \eqref{eq:oneillspin} and
\eqref{eq:gauss} that $A(X)=X$ for all $X$ orthogonal to $\xi.$
Thus, the manifold $N$ is flat and by the de Rham theorem the
universal cover $\widetilde{M}$ of $M$ is isometric to the
Riemannian product $\mathbb{R}\times Z,$ where $Z$ is a simply
connected Riemannian manifold carrying a maximal number of Killing
spinors. Therefore with the use of B\"ar's classification
\cite{B2}, this implies in particular that $\widetilde{M}\simeq
\mathbb{R}\times \mathbb{S}^n.$ Moreover, the fundamental group of
$M$ is embedded in the product ${\rm Isom}_+(\R)\times {\rm
Isom}_+(\mathbb{S}^n)$ where ${\rm Isom}_+$ denotes the group of
isometries which preserves the orientation of the corresponding
manifold.
Since $n$ is even, we can deduce that $\pi_1(M)\cong\mathbb{Z}$ and acts on $\mathbb{R}\times\mathbb{S}^n$ via $\left(k,(t,x)\right)\mapsto\left(t+ka,A^kx\right)$ for some $(a,A)\in\mathbb{R}^\times\times\mathrm{SO}_{n+1}$.
This action lifts to the spin level and induces two spin structures on $M$ which both admit a maximal number of linearly independent $-\frac{1}{2}$- and $\frac{1}{2}$-basic Killing spinors, see e.g. \cite[Thm. 3.2]{AlexGrantchIv98}.

\noindent In order to check that $N$ is isometric to $\lquot{\mathbb{R}\times B}{\mathbb{Z}}$, where $B$ is the closed unit ball, we need the following lemma:
\begin{lemma}\label{extension}
Under the same conditions as in {\rm Theorem \ref{thm:6}}, the vector field $\xi$ can be extended to a unique parallel vector field $\hat \xi$ on $N$.
In particular, $\hat{\xi}$ is orthogonal to the unit normal $\nu$ along the boundary $M$.
\end{lemma}
Assume this lemma for a moment and consider a connected integral
submanifold $N_1$ of the bundle $(\R\hat{\xi})^\perp,$ where the
orthogonal is taken in $N$. From the parallelism of the vector
field $\hat{\xi}$, it is straightforward to see that $N_1$ is the
quotient of a totally geodesic hypersurface $\widetilde{N}_1$ of
the universal cover $\widetilde{N}$ (which is complete) of $N$. In
particular, the manifold $\widetilde{N}_1$ is complete since it is
a level hypersurface of the function $f$ defined on
$\widetilde{N}$ by $d^{\widetilde{N}} f=\widetilde{\hat\xi}$
(recall here that $d^N\hat\xi=0$). From the fact that the
universal cover is a local isometry, we deduce that $N_1$ is
complete. On the other hand, the boundary of $N_1$ is also a
totally geodesic hypersurface in $\partial N=\lquot{\mathbb{R}\times\mathbb{S}^n}{\mathbb{Z}}$ with normal vector field $\xi$,
carrying a maximal number of Killing spinors (hence it is compact
but may be unconnected). Moreover, the second fundamental form of
$\partial N_1$ is equal to $-\nabla_X^{N_1}\nu=-\nabla_X^N\nu=X,$
which means that $\partial N_1$ is totally umbilical in $N_1$. Thus
from the rigidity result in \cite[Thm. 1.1]{Li}, we deduce that
$N_1$ is compact which from \cite[Cor. 4]{HM1} or \cite[Thm.
1]{Ra} implies that $\partial N_1$ is connected and is isometric
to the round sphere. Therefore $N_1$ -- and thus $\widetilde{N}_1$ -- is isometric to the unit ball
$B$ (see \cite[Cor. 5]{HM1} or \cite[Thm. 5]{Ra}).
Finally, by the de Rham theorem, the manifold $\widetilde{N}$ is isometric to $\mathbb{R}\times B$ and therefore $N$ is the quotient of
the Riemannian product $\R\times B$ by its fundamental group.
Since $\pi_1(N)$ embeds into $\pi_1(M)$ (any isometry of $B$ fixing pointwise $\mathbb{S}^n$ is the identity map), $N$ is isometric to $\lquot{\mathbb{R} \times B}{\mathbb{Z}}$.
\hfill$\square$\\ \\
Now, let us prove Lemma \ref{extension} that we used to prove Theorem \ref{thm:6} for $n$ even. \\ \\
{\bf Proof of Lemma \ref{extension}.} Given any $p$-form $\omega$ on $M$ the following boundary problem
\begin{equation} \label{eq:boundaryproblem}
\left\{\begin{array}{ll}
\Delta^N \hat{\omega}=0 &\textrm{on $N$},\\\\
J^*\hat\omega=\omega,\,\, J^*(\delta^N\hat\omega)=0 &\textrm{on $M$}
\end{array}\right.
\end{equation}
admits a solution $\hat\omega$ (which is also a $p$-form) on $N$
(see Lemma 3.5.6 in \cite{S}). Here $J^* \hat\omega$ denotes the
restriction of $\hat\omega$ to the boundary. Using Stokes formula,
on can easily check that $\delta^N d^N\hat\omega=0$ (see e.g. the
proof of \cite[Theorem 8]{RS}). Therefore, the form
$\hat\phi:=d^N\hat\omega$ satisfies the boundary problem
$$
\left\{\begin{array}{ll}
d^N \hat{\phi}=\delta^N\hat\phi=0 &\textrm{on $N$}, \\\\
J^*\hat\phi=d^M\omega &\textrm{on $M.$}
\end{array}\right.
$$

\noindent Given the $1$-form $\omega=\xi^\flat$ on $M$, let $\hat\xi$ be a solution of the problem \eqref{eq:boundaryproblem}. In the sequel, we will show that the vector field $\hat\xi$ is parallel on the manifold $N$. To do this, we will use the Reilly formula established in \cite[Theorem 3]{RS1}: For any $p$-form $\alpha$ on $N$ with $p\geq 1$, we have
$$\int_N|d^N\alpha|^2+|\delta^N\alpha|^2=\int_N|\nabla^N\alpha|^2+\langle W^{[p]}(\alpha),\alpha\rangle+2\int_M\langle\nu\lrcorner\alpha,\delta^M(J^*\alpha)\rangle+\int_M \mathcal{T}(\alpha,\alpha)$$
where $W^{[p]}$ is the curvature term of $N$ in the Bochner-Weitzenb\"ock formula (which in our case vanishes since $N$ is flat) and
$$\mathcal{T}(\alpha,\alpha)=\langle A^{[p]}(J^*\alpha),J^*\alpha\rangle+(n+1)H|\nu\lrcorner\alpha|^2-\langle A^{[p-1]}(\nu\lrcorner\alpha),\nu\lrcorner\alpha\rangle,$$
with $(A^{[k]}\beta )(X_1,\dots,X_k)=\sum_{i=1}^k \beta(X_1,\dots, A(X_i),\dots,X_k)$ where $A$ is the second fundamental form (here $\beta$ is a  $k$-form on $N$). By convention, we take $A^{[0]}=0$.

\noindent The form $\hat\phi=d^N\hat\xi$ is closed and co-closed on $N$ and its restriction to the boundary vanishes from the fact that $d^M\xi=0$ (the vector field $\xi$ is parallel on $M$). Therefore, the Reilly formula  applied to $\hat\phi=d^N\hat\xi$ gives that
$$0=\int_N|\nabla^N\hat\phi|^2+(n+1)\int_M H|\nu\lrcorner\hat\phi|^2-\int_M\langle A^{[1]}(\nu\lrcorner\hat\phi),\nu\lrcorner\hat\phi\rangle.$$
Since the second fundamental form is equal to zero in the direction of $\xi$ and to the identity in the orthogonal direction to $\xi$, we find for $A^{[1]}$ that
$$(A^{[1]}\beta) (X)=\beta(X)\,\, \text{for}\,\, X\perp \xi\,\, \text{and}\,\, (A^{[1]}\beta) (\xi)=0.$$


\noindent Hence, we get that
$$\langle A^{[1]}(\nu\lrcorner\hat\phi),\nu\lrcorner\hat\phi\rangle=\sum_{i=1}^n(\nu\lrcorner\hat\phi)(e_i)^2,$$
where $\{e_i\}_{i=1,\cdots,n}$ is an orthonormal frame of $\Gamma(Q)$. Thus, after replacing the mean curvature $H$ by $\frac{n}{n+1}$, the Reilly formula reduces to 
$$0=\int_N |\nabla^N\hat\phi|^2+(n-1)\int_M |\nu\lrcorner \hat\phi|^2+\int_M(\nu\lrcorner \hat\phi)(\xi)^2.$$ 

Therefore the form $\nu\lrcorner \hat\phi$ vanishes on $M$ and $\hat\phi$ is parallel on $N$, which means that it has a constant norm. But using the identity $|\hat\phi|^2=|J^*\hat\phi|^2+|\nu\lrcorner \hat\phi|^2$ which holds at any point on the boundary, one can deduce that $\hat\phi=0$ (recall here that the restriction of $\hat\phi$ to the boundary vanishes). Thus we get that $d^N\hat\xi=0$. On the other hand, with the help of the Stokes formula and the fact that $\hat\xi$ is a solution of the problem \eqref{eq:boundaryproblem}, one can easily prove that $\delta^N\hat\xi=0$. We apply again the Reilly formula for the $1$-form $\hat\xi$ to get that 
$$0=\int_N|\nabla^N\hat\xi|^2+n\int_M|\nu\lrcorner \hat\xi|^2,$$ 
which gives that $\hat\xi$ is parallel on $N$ and is orthogonal to $\nu$ along the boundary. Here we used that $\delta^M\xi=0$ and $A^{[1]}(\xi)=0.$ 

\noindent Finally, we would like to notice that the vector field $\hat\xi$ is unique. Indeed, we prove that any extension of the vector field $\xi$ of the boundary problem \eqref{eq:boundaryproblem} is a parallel vector field on the whole manifold. Hence for any other extension $\hat\xi_1,$ the vector field $\hat\xi-\hat\xi_1$ has a constant norm (both are parallel) and restricts to zero on the boundary.
\hfill$\square$

\noindent{\bf Proof of Theorem \ref{thm:6} for $n$ odd.} The existence of a maximal number of basic Killing spinors of constant $-\frac{1}{2}$ gives $2^{[\frac{n}{2}]}$ parallel spinors on $N,$ by Theorem \ref{thm:5}. However, given a basic Killing spinor $\varphi$ of constant $\frac{1}{2}$ which belongs to $\Sigma Q\simeq \Sigma^+M$ implies by the identity \eqref{eq:commute} that $\xi\cdot_M\varphi$ is a solution of the basic Dirac equation belonging to $\Sigma Q\simeq \Sigma^-M.$ Therefore, Theorem \ref{thm:5} assures again the existence of $2^{[\frac{n}{2}]}$ parallel spinors on $N$ which with the first family provides a maximal number of parallel spinors. This forces $N$ to be flat. The proof carries over as in the even case with two differences: first, the only compact, connected manifold carrying a maximal number of Killing spinors for both constants $-\frac{1}{2}$ and $\frac{1}{2}$ is the round sphere \cite[Thm. 4]{B3}; second, one cannot deduce that $\Gamma=\pi_1(M)$ is isomorphic to $\mathbb{Z}$.
\hfill$\square$

\begin{remarks}
{\rm\noindent\begin{enumerate}
\item We notice that the Killing constant $-\frac{1}{2}$ in Theorem \ref{thm:6} can be replaced by any constant $-\alpha$. In this case, the term $\frac{n}{n+1}$ in the inequality $\frac{n}{n+1} + \frac{1}{n+1} [\frac n2]^{\frac 12} \vert \Omega\vert \leq H$ should be replaced by $\frac{2n\alpha}{n+1}$. On the other hand, the condition taken on the mean curvature cannot be dropped. In fact, consider in dimension $4$ the example of the unit closed ball where the boundary $\mathbb{S}^3$ is endowed with the unit Killing vector field which defines the Hopf fibration over the sphere $\mathbb{S}^2$ of radius $\frac{1}{2}.$ The normal bundle of the flow (which is isometric to the tangent bundle of $\mathbb{S}^2$) carries a K\"ahler structure and has a $2$-dimensional space of basic Killing spinors of constants $\pm 1$. However, we have $\frac{2n\alpha}{n+1}+\frac{1}{n+1}[\frac{n}{2}]^\frac{1}{2}|\Omega|=\frac{4}{3}+\frac{1}{3}(1)=\frac{5}{3}>H=1.$
\item In case $n$ is odd, the subgroup $\Gamma$ by which $\mathbb{R}\times\mathbb{S}^n$ is modded out is not necessarily isomorphic to $\mathbb{Z}$.
Consider for instance $\Gamma:=\mathbb{Z}\times\Gamma_2$, where $\Gamma_2\subset\mathrm{SO}_{n+1}$ is a finite subgroup consisting of rotations in orthogonal $2$-planes in $\mathbb{R}^{n+1}$. 
Then $\Gamma$ acts freely and properly discontinuously on $\mathbb{R}\times\mathbb{S}^n$ with compact quotient.
\end{enumerate}
}
\end{remarks}

\noindent An immediate consequence of Theorem \ref{thm:6} is when
the boundary of a spin manifold is isometric to the Riemannian
product $\mathbb{S}^1\times \mathbb{S}^n$ (which is a trivial
fibration over the round sphere). In this case, if the spin
structure on $N$ induces  the trivial one on the boundary
$\mathbb{S}^1\times \mathbb{S}^n$, the flow carries obviously a
maximal number of basic Killing spinors and Theorem \ref{thm:6} applies.
In fact, we can say more:

\begin{cor}\label{c:S1B} Let $(N^{n+2},g)$ be a compact spin Riemannian manifold with non-negative scalar curvature.
We assume that the boundary is isometric to $\mathbb{S}^1\times
\mathbb{S}^n$ with mean curvature $H\geq \frac{n}{n+1}$. If the
induced spin structure on $M$ is the trivial one on
$\mathbb{S}^1\times \mathbb{S}^n$, then $N$ is isometric to the
product of $\mathbb{S}^1$ with the unit ball.
\end{cor}

{\bf Proof.} As a consequence of Theorem \ref{thm:6}, $M$ must be isometric to the quotient of $\mathbb{R}\times\mathbb{S}^n$ by some discrete fixed-point-free cocompact subgroup $\Gamma\subset\mathbb{R}\times\mathrm{SO}_{n+1}$. 
But since $M$ is \emph{by assumption} isometric to $\mathbb{S}^1\times\mathbb{S}^n$, the $\Gamma$-action must be trivial on the $\mathbb{S}^n$-factor, so that $N$ is isometric to $\mathbb{S}^1\times B$.
\hfill$\square$

\subsection{Integrability condition of the normal bundle}

\noindent In the following, we will state another rigidity results which mainly uses a previous work of B.-Y.~Chen \cite{Ch}. In his paper, Chen proved that given any Riemannian submersion $M^{n+1}\rightarrow B^b$ with totally geodesic fibres where $M$ is isometrically immersed onto a Riemannian manifold $N$, the O'Neill tensor of the submersion can be bounded from above by the mean curvature of the immersion and the sectional curvature of $N$. Indeed, he showed that
$$|h|^2\leq \frac{(n+1)^2}{4}H^2+b(n+1-b)\ {\rm max}\, \tilde{K}(p),$$
where ${\rm max}\, \tilde{K}(p)$ denotes the maximum value of the
sectional curvature of $N$ restricted to plane sections in $T_pM.$
Using this estimate, we can state the following:

\begin{thm}\label{final}
Let $(N,g)$ be an $(n+2)$-Riemannian spin manifold of nonnegative
scalar curvature with connected boundary $M$ of mean curvature
satisfying
$H\geq\frac{2n\sqrt{2}}{(n+1)(2\sqrt{2}-[\frac{n}{2}]^\frac{1}{2})}$
for $n <16$ such that the sectional curvature vanishes on the
boundary $M$. Assume that $M$ is endowed with a Riemannian flow
given by a unit vector field $\xi.$ If there exists a solution
$\varphi$ of $D_b\varphi=\frac{n}{2}\varphi,$ the normal bundle
$Q$ is integrable.

\end{thm}

\noindent{\bf Proof.} Suppose that the normal bundle is not integrable, that is $h$ is non-zero. By Chen's estimate we get
$$|h|=\sqrt{2}|\Omega|\leq \frac{n+1}{2}H\leq\frac{\sqrt{2}}{[\frac{n}{2}]^\frac{1}{2}}((n+1)H-n),$$
which gives a contradiction from Corollary \ref{thm:rig}.
\hfill$\square$

{\bf Acknowledgment.} The authors are  indebted to Oussama Hijazi
and Aziz El Kacimi for many enlightening discussions and comments.
The fourth named author gratefully acknowledges the financial
support of the Berlin Mathematical School (BMS) and  would like to
thank the University of Potsdam, especially Christian B\"ar and
his group, for their generous support and friendly welcome.
We also thank Christian B\"ar for pointing out a mistake in the proof of Theorem \ref{thm:6}.

\end{document}